\documentclass[12pt]{amsart}
\usepackage{amssymb}
\usepackage[all]{xy}

 \rm

\renewcommand{\(}{\left(}
\renewcommand{\)}{\right)}

\renewcommand{\>}{\rangle}
\newcommand{\x}{\times}
\renewcommand{\bar}{\overline}
\newcommand{\abs}[1]{\left\lvert#1\right\rvert}

\newcommand{\st}{\:|\:}

\newcommand{\C}{{\mathbb{C}}}

\renewcommand{\phi}{\varphi}

\renewcommand{\H}{{\mathcal{H}}}
\newcommand{\BH}{{\mathcal{B}}(\H)}
\newcommand{\A}{{\mathcal{A}}}
\newcommand{\B}{{\mathcal{B}}}

\theoremstyle{plain}
\newtheorem{thm}{Theorem}
\newtheorem{lem}[thm]{Lemma}

\newtheorem{cor}[thm]{Corollary}

\theoremstyle{definition}

\theoremstyle{remark}

\title{Inductive algebras for the affine group of a finite field}

\author{Promod~Sharma}
\author{M.~K.~Vemuri}
\address{Department of Mathematical Sciences\\ IIT (BHU)\\ Varanasi 221 005\\
INDIA}

\begin{document}

\begin{abstract}
Each irreducible representation of the affine group of a finite
field has a unique maximal inductive algebra, and it is self adjoint.
\end{abstract}


\maketitle
\thispagestyle{empty}


\section{Introduction}\label{S:intro}
Let $G$ be a separable locally compact group and $\pi$ an irreducible
unitary representation of $G$ on a separable Hilbert space $\H$.  Let
$\BH$ denote the algebra of bounded operators on $\H$.  An {\em
inductive algebra} is a weakly closed abelian sub-algebra $\A$ of $\BH$
that is normalized by $\pi(G)$, i.e., $\pi(g) \A \pi(g)^{-1} = \A$ for each
$g \in G$.  If we wish to emphasize the
dependence on $\pi$, we will use the term $\pi$-inductive algebra.
A {\em maximal inductive algebra} is a maximal element of the set of
inductive algebras, partially ordered by inclusion.

The identification of inductive algebras can shed light on the
possible realizations of $\H$ as a space of sections of a homogeneous
vector bundle (see e.g. \cite{Stegel-2004,Stegel-2006,sl2r,rcr}).
For self-adjoint maximal inductive algebras there is a precise result
known as Mackey's Imprimitivity Theorem, as explained in the
introduction to \cite{Stegel-2004}.  Inductive algebras have also found
applications in operator theory (see e.g. \cite{iahs, Koranyi-2014}).

In \cite{Raghavan-2005}, it was shown that finite dimensional
inductive algebras for a connected group are trivial.  However, there
are interesting inductive algebras for finite groups (see e.g.
\cite{iafhg}).  In this note, we classify the maximal inductive algebras
for the representations of the affine group (the ``$ax+b$'' group)
of a finite field.

In Section \ref{S:ax+b}, we recall the structure of the affine group
of a finite field, and set up the notation.  In Section
\ref{S:representations}, we recall its representation theory, and
formulate our main result.  The main result is proved in
Section \ref{S:proof}.

\section{The affine group of a finite field}\label{S:ax+b}
Let $k$ be a finite field of order $q=p^n$, where $p$ is prime.
Let $k^\times$ denote the
multiplicative group of non-zero elements of $k$.  Recall that the
affine group of $k$ is the group $G$ of affine-automorphisms of $k$.
Thus an element $g$ of $G$ is a map $g: k \to k$ of the form $g(x) =
ax+b$ where $a \in k^\times$ and $b \in k$, and the group law is
composition.  The group $G$ may be identified with the group of matrices
\begin{equation*}
\left\{
\left.
\left[
\begin{matrix}
a & b \\
0 & 1
\end{matrix}
\right] 
\; \right | \:
a \in k^\times, \, b \in k
\right\}.
\end{equation*}

Let $\iota: k \to G$, $p: G \to k^\times$ and $s: k^\times \to G$  be defined by
\begin{equation*}
\iota(b) = 
\left[
\begin{matrix}
1 & b \\
0 & 1
\end{matrix}
\right], \quad
p\(
\left[
\begin{matrix}
a & b \\
0 & 1
\end{matrix}
\right] 
\) =
a, \;\text{and}\quad
s(a) =
\left[
\begin{matrix}
a & 0 \\
0 & 1
\end{matrix}
\right] .
\end{equation*}
Then $\iota, p$ and $s$ are homomorphisms, $\iota(k) \triangleleft G$ and
\begin{equation*}
\xymatrix{
0 \ar[r] & k \ar[r]^\iota & G \ar[r]^p
                              \ar@/_/@{<--}[r]_s  & k^\times \ar[r] & 1
}
\end{equation*}
is an exact sequence with splitting $s$.  Thus $G$ is a semidirect
product $k^\x \ltimes k$.  We note for future reference that
$s(a')\iota(b)s(a')^{-1} = \iota(a'b)$.

\medskip

\section{The representations and their inductive algebras}
\label{S:representations}

The irreducible unitary representations of $G$ may be constructed using
the {\em Mackey machine} (see \cite[\S 3.9]{Mackey-1976}).  There are
$q-1$ characters (one dimensional representations), and one
$(q-1)$-dimensional representation (up to unitary equivalence).

Obviously, the characters have only the trivial
inductive algebra $\C$, which is self adjoint.

The $(q-1)$-dimensional representation is
\begin{equation*}
\pi = \mathrm{Ind}_{\iota(k)}^G \chi,
\end{equation*}
where $\chi: k \to \C^\times$ is a non-trivial homomorphism
(i.e. $\chi \not\equiv 1$).
Let $\H$ denote the Hilbert space of all complex valued
functions on $k^\times$ equipped with the inner product
\begin{equation*}
\<F_1,F_2\> = \sum_{a'\in k^\times} F_1(a')\bar{F_2(a')}.
\end{equation*}
The representation $\pi$ may be realized on $\H$ by
\begin{equation*}\label{G-action}
(\pi(g) F)(a')
=  \chi(a'b) F(a'a), \qquad g=\iota(b)s(a).\\
\end{equation*}

For each $\phi \in \H$, let $m_\phi: \H\to\H$ be
defined by $m_\phi(F)=\phi F$.  Let
\begin{equation*}
\B=\{m_\phi \st \phi \in \H\}.
\end{equation*}
Then $\B$ is a maximal-abelian subalgebra of 
$\BH$, and that $\B$ is $\pi$-inductive.  Therefore $\B$ is a maximal
$\pi$-inductive algebra.  Moreover, it is self-adjoint.
Our main result is the following theorem.

\begin{thm}
$\B$ is the only maximal $\pi$-inductive algebra.
\end{thm}

\section{The proof}\label{S:proof}

Let $\A$ be a maximal $\pi$-inductive algebra.

\begin{lem}\label{L:reduced}
There are no non-zero nilpotent elements in $\A$.
\end{lem}

\begin{proof}
Let $\mathcal{N}$ denote the set of nilpotent elements in $\A$ (the
nil-radical of $\A$).  Let
\begin{equation*}
\mathcal{K}=\{F \in \H \st T(F)=0, \quad \forall T \in \mathcal{N}\}.
\end{equation*}
By (a trivial case of) Engel's theorem
\cite[\S3.3]{Humphreys-1972}, $\mathcal{K}\ne 0$.  Observe that
$\mathcal{N}$ is normalized by $\pi(G)$, so $\mathcal{K}$ is
$\pi(G)$-invariant.  However, since $\pi$ is irreducible, it follows
that $\mathcal{K}=\H$, whence $\mathcal{N}=0$.
\end{proof}

\begin{cor}\label{C:dimension-bound}
$\mathrm{dim}(\A) \le q-1.$
\end{cor}

\begin{proof}
By Lemma \ref{L:reduced}, the Jordan-Chevalley decomposition
\cite[\S4.2]{Humphreys-1972}, and
the fact that $\A$ is abelian, it follows that there is a
(not necessarily orthonormal) basis for
$\H$ in which each element of $\A$ is diagonal.  Since
$\mathrm{dim}(\H) = q-1$, the result follows.
\end{proof}

For $b' \in k$, define $\kappa(b'):\A\to\A$ by
$\kappa(b')T = \pi(\iota(b'))T\pi(\iota(b'))^{-1}$.  Then $\kappa$ is
a representation of the finite abelian group $k$ on the vector space $\A$,
which decomposes as
\begin{equation*}
\A = \bigoplus_{b \in k} \A_b.
\end{equation*}
where
\begin{equation*}
\A_b = \{T \in \A \st \kappa(b')T = \chi(bb') T,
                      \quad\forall b' \in k\}.
\end{equation*}
Here we are using the fact that every character of $k$ is of the form
$\chi_b$ where $\chi_b(b')=\chi(bb')$ for all $b'\in k$.

Observe that
\begin{enumerate}
\item if $T \in \A_b$ and $T' \in \A_{b'}$, then $TT' \in \A_{b+b'}$, and
\item for each $a\in k^\times$, the map
$T\mapsto \pi(s(a))^{-1} T \pi(s(a))$ is a linear isomorphism
$\A_b \to \A_{ab}$.
\end{enumerate}

\begin{lem}\label{L:A1}
$\A_1=0.$
\end{lem}

\begin{proof}
Suppose not.  Then there exists a non-zero element $T \in \A_1$.
By the first observation above,
$T^p \in \A_{0}$.  Since $T$ is not nilpotent (see Lemma \ref{L:reduced}),
it follows that $\mathrm{dim}(\A_0) \ge 1$.  For $b\ne 0$, we have
$\mathrm{dim}(\A_b) = \mathrm{dim}(\A_1) \ge 1$, by the second observation.
Therefore
\begin{equation*}
\dim(\A) = \sum_{b\in k} \mathrm{dim}(\A_b) \ge \abs{k}=q,
\end{equation*}
contradicting Corollary \ref{C:dimension-bound}.
\end{proof}

\begin{lem}
$\A_0 \subseteq \B$.
\end{lem}

\begin{proof}
Observe that $\H$ is spanned by the set $\{\chi_{b'}|_{k^\times} \st b' \in k\}$,
hence
$\B$ is spanned by the set $\{m_{(\chi_{b'}|_{k^\times})} \st b' \in k\}$.  If
$T \in \A_0$,
then $T$ commutes with $m_{(\chi_{b'}|_{k^\times})}$ for each $b'\in k$, whence $T$
commutes with $\B$.  Since $\B$ is maximal-abelian, it follows that $T\in\B$.
\end{proof}

By Lemma \ref{L:A1}, and the second observation above, it follows that
$\A_b=0$ for all $b \in k^\times$, whence $\A=\A_0 \subseteq \B$.
Since $\A$ is maximal, it follows that $\A=\B$.

\bibliographystyle{amsplain}
\bibliography{v4-iaftagoaff.bib}

\end{document}